\magnification=\magstephalf
\baselineskip=15pt

\def\qed{\hbox{\vrule height5pt width5pt}}
\parindent 0in

\bigskip
\centerline{{\bf Conditionally Positive Functions and p-norm Distance Matrices }}
\medskip
\centerline{B. J. C. Baxter\footnote{${}^*$}{DAMTP,
Silver Street, Cambridge
CB3 9EW, England.}}
\smallskip
\vskip 15mm \noindent

{\bf Abstract.}\quad In Micchelli [3], deep results were obtained concerning the
invertibility of matrices arising from radial basis function
interpolation. In particular, the Euclidean distance matrix
was shown to be invertible for distinct data. In this paper, we investigate the invertibility
of distance matrices generated by $p$-norms.  In particular, we show
that, for any $p\in (1, 2)$, and for distinct points $ x^1,
..., x^n \in {\cal R}^d $, where $n$ and $d$ may be any positive
integers, with the proviso that $ n \ge 2$, the matrix $A \in
{\cal R}^{n \times n} $ defined by $$ A_{ij} = \Vert x^i - x^j \Vert_p ,
\hbox{ for } 1 \le i, j \le n, $$
 satisfies $$ (-1)^{n-1}\det A > 0 .$$

We also show how to construct, for every $p > 2$, a
configuration of distinct points in some ${\cal R}^d$ giving a
singular $p$-norm distance matrix. Thus radial basis function interpolation
using $p$-norms is uniquely determined by any distinct data for
$p \in (1,2]$, but not so for $p > 2$.

\vskip 36pt
\leftline {\sevenrm 1980 AMS classification: Primary 41A05, 41A63, 41A25; Secondary 41A30.}
\leftline {\sevenrm Keywords and phrases: Multivariate Interpolation, Multivariate Approximation, Radial Basis Functions, p-norms.}

\vfill
\eject

\def\qed{\hbox{\vrule height5pt width5pt}}
\parindent 0in

%**end of header

{\bf Introduction }

\vskip 6pt

\hskip 12pt The real multivariate interpolation problem is as follows. Given distinct points
$x^1, \ldots , x^n \in {\cal R}^d$ and real scalars $f_1, \ldots ,
f_n$, we wish to construct a continuous function $s:{\cal R}^d
\rightarrow {\cal R}$ for which $$s(x^i) = f_i, \hbox{ for } i = 1,
\ldots ,n .$$ The radial basis function approach is to choose a
function $\varphi :[0, \infty ) \rightarrow [0, \infty )$ and a norm
$\Vert . \Vert $ on ${\cal R}^d$ and then let $s$ take the form $$
s(x) = \sum_{i = 1}^n \lambda_i \/ \varphi (\Vert x - x^i \Vert) .$$
Thus $s$ is chosen to be an element of the vector space spanned by the
functions $\xi \mapsto \varphi (\Vert \xi - x^i \Vert )$, for $i = 1,
\ldots , n $.  The interpolation conditions then define a linear
system $ A \lambda = f$, where $A \in {\cal R}^{n \times n}$ is given
by $$ A_{ij} = \varphi ( \Vert x^i - x^j \Vert ), \hbox{ for } 1 \le
i, j \le n,$$ and where $\lambda = ( \lambda_1, ..., \lambda_n )$ and
$f = (f_1, ...,f_n)$. In this paper, a matrix such as $A$ will be
called a distance matrix.

\hskip 12pt Usually $\Vert . \Vert $ is chosen to be the Euclidean norm,
and in this case Micchelli [4] has shown the distance matrix generated
by distinct points to be invertible for several useful choices of
$\varphi $.  In this paper, we investigate the invertibility of the
distance matrix when $\Vert . \Vert $ is a $p$-norm for $1 < p <
\infty$, $p
\ne 2$, and $\varphi (t) = t$, the identity. We find that $p$-norms do
indeed provide invertible distance matrices given distinct points, for
$1 < p \le 2$. Of course, $p = 2$ is the Euclidean case mentioned
above and is not included here.  Now Dyn, Light and Cheney [2] have
shown that the $1-$norm distance matrix may be singular on quite
innocuous sets of distinct points, so that it might be useful to
approximate $\Vert . \Vert_1$ by $\Vert . \Vert_p$ for some $p \in (1,
2]$.  This work comprises section 2. The framework of the proof is
very much that of Micchelli [4].

\hskip 12pt For every $p > 2$, we find that distance matrices can be singular
on certain sets of distinct points, which we construct. We find that
the higher the dimension of the underlying vector space for the points
$x^1, \ldots , x^n$, the smaller the least $p$ for which there exists
a singular $p$-norm.

\vskip 18pt

{\bf 1. Almost negative matrices}

\vskip 12pt

\hskip 18pt Almost every matrix considered in this paper will induce a
non-positive form on a certain hyperplane in ${\cal R}^n$.
Accordingly, we first define this ubiquitous subspace and fix
notation.

\proclaim {Definition 1.1}. For any positive integer $n$, let
$$ Z_n = \{ \ y \in {\cal R}^n : \sum_{i=1}^n y_i = 0 \ \}. $$

Thus $Z_n$ is a hyperplane in ${\cal R}^n $. We note that $Z_1 = \{ 0
\}$.

\vskip 18pt

\proclaim {Definition 1.2}. We shall call $A \in {\cal R}^{n \times n} $ {\bf almost negative definite} (AND) 
if $A$ is symmetric and $$ y^T Ay \le 0, \hbox{ whenever } y \in
Z_n.$$ Furthermore, if this inequality is strict for all non-zero $y
\in Z_n$, then we shall call $A$ {\bf strictly AND}.
\vskip 18pt

\proclaim {Proposition 1.3}. Let $A \in {\cal R}^{n \times n} $ be
strictly AND with non-negative trace. Then $$ (-1)^{n-1} \det A > 0.$$

{\it Proof. }We remark that there are no strictly AND $1 \times 1$
matrices, and hence $n \ge 2$. Thus $A$ is a symmetric matrix inducing
a negative-definite form on a subspace of dimension $n-1 > 0$, so that
$A$ has at least $n-1$ negative eigenvalues. But trace $A \ge 0$, and
the remaining eigenvalue must therefore be positive \qed
\vskip 18pt
Micchelli [4] has shown that both $A_{ij} = |x^i-x^j|$ and $A_{ij} =
(1+|x^i-x^j|^2)^{1\over 2}$ are AND, where here and subsequently $|.|$
denotes the Euclidean norm.  In fact, if the points $x^1, \ldots, x^n$
are distinct and $n \ge 2$, then these matrices are strictly AND.
Thus the Euclidean and multiquadric interpolation matrices generated
by distinct points satisfy the conditions for proposition 1.3.
\vskip 12pt

\hskip 18pt The work in this paper rests on the following characterization of
AND matrices with all diagonal entries zero. This theorem is stated
and used to good effect by Micchelli [4], who omits much of the proof
and refers us to Schoenberg [5]. Because of its extensive use in this
paper, we include a proof for the convenience of the reader.  The
derivation follows the same lines as that of Schoenberg [5].

\vskip 18pt

\proclaim {Theorem 1.4}. Let $A \in {\cal R}^{n \times n} $ have all diagonal entries zero. Then $A$ is
AND if and only if there exist $n$ vectors $ y^1, \ldots, y^n \in
{\cal R}^n $ for which $$ A_{ij} = |y^i - y^j|^2. $$

{\it Proof.} Suppose $ A_{ij} = |y^i - y^j|^2 $ for vectors $ y^1,
\ldots, y^n \in {\cal R}^n$. Then $A$ is symmetric and the
following calculation completes the proof that $A$ is AND. Given any $z
\in Z_n$, we have $$\eqalign { z^T A z &= \sum_{i,j=1}^n z_i z_j |y^i -
y^j|^2 \cr &= \sum_{i,j=1}^n z_i z_j ( |y^i|^2 + |y^j|^2 -
2(y^i)^T(y^j) )\cr &= -2 \sum_{i,j=1}^n z_i z_j (y^i)^T(y^j) \ ,
\hbox{ since the coordinates of $z$ sum to zero, } \cr &= -2 \ |
\sum_{i=1}^n z_i y^i \ |^2 \le 0 .}$$

This part of the proof is given in Micchelli [4].  The converse
requires two lemmata.

\proclaim {Lemma 1.5}. Let $B \in {\cal R}^{k \times k} $ be a symmetric non-negative definite matrix. 
Then we can find $\xi^1,
\ldots, \xi^k \in {\cal R}^k $ such that $$ B_{ij} = |\xi^i|^2 +
|\xi^j|^2 - |\xi^i - \xi^j|^2 . $$

{\it Proof.} Since $B$ is symmetric and non-negative definite, we have
$ B = P^T P$, for some $P\in {\cal R}^{k \times k} $.

Let $p^1, \ldots, p^k$ be the columns of $P$. Thus $$ B_{ij} = (p^i)^T
(p^j).$$

Now $$ |p^i - p^j|^2 = |p^i|^2 + |p^j|^2 - 2(p^i)^T(p^j).$$

Hence $$ B_{ij} = {1 \over 2}( |p^i|^2 + |p^j|^2 - |p^i - p^j|^2 ).$$

All that remains is to define $ \xi^i = p^i / \surd 2$ , for $i = 1,
\ldots, k$ \qed

\vskip 18pt

\proclaim {Lemma 1.6}. Let $A \in {\cal R}^{n \times n} $. Let $e^1, \ldots, e^n$ denote
the standard basis for ${\cal R}^n $, and define $$\eqalign { f^i &=
e^n - e^i ,\hbox{ for } i = 1, \ldots, n-1, \cr f^n &= e^n.}$$
Finally, let $F \in {\cal R}^{n \times n}$ be the matrix with columns
$ f^1, \ldots, f^n.$ Then $$\eqalign { (-F^T AF)_{ij} &= A_{in} +
A_{nj} - A_{ij} - A_{nn}, \hbox{ for } 1 \le i, j \le n-1 \ , \cr
(-F^T AF)_{in} &= A_{in} - A_{nn}, \cr (-F^T AF)_{ni} &= A_{ni} -
A_{nn}, \hbox{ for } 1 \le i \le n-1 \ ,\cr (-F^T AF)_{nn} &= -A_{nn}.
}$$

{\it Proof.} We simply calculate $(-F^T AF)_{ij} \equiv -(f^i)^T A
(f^j)$ \qed
\vskip 18pt

\hskip 18pt We now return to the proof of theorem 1.4: Let $A \in {\cal R}^{n
\times n} $ be AND with all diagonal entries zero.  Lemma
1.6 provides a convenient basis from which to view the action of $A$.
Indeed, if we set $ B = -F^T AF$, as in lemma 1.6, we see that the
principal submatrix of order $n-1$ is non-negative definite, since
$f^1, \ldots, f^{n-1}$ form a basis for $Z_n$. Now we appeal to Lemma
1.5, obtaining $ \xi^1,\ldots, \xi^{n-1} \in {\cal R}^{n-1}$ such that
$$ B_{ij} = |\xi^i|^2 + |\xi^j|^2 - |\xi^i -
\xi^j|^2 \ , \hbox{ for } 1 \le i, j \le n-1, $$

while lemma 1.6 gives $$ B_{ij} = A_{in} + A_{jn} - A_{ij}.$$

Setting $i = j$ and recalling that $A_{ii} = 0$, we find $$ A_{in} =
|\xi^i|^2 , \hbox{ \ for } 1 \le i\le n-1$$ and thus we obtain $$
A_{ij} = |\xi^i - \xi^j|^2, \hbox{ for } 1 \le i, j \le n-1 .$$

Now define $\xi^n = 0$. Thus $ A_{ij} = |\xi^i - \xi^j|^2$, for $1 \le
i, j \le n$, where $ \xi^1,\ldots, \xi^n \in {\cal R}^{n-1}$. We may
of course embed ${\cal R}^{n-1}$ in ${\cal R}^{n}$. More formally, let
$\iota : {\cal R}^{n-1} \hookrightarrow {\cal R}^{n} $ be the map
$\iota:(x_1, \ldots, x_{n-1}) \mapsto (x_1, \ldots, x_{n-1}, 0)$, and,
for $i = 1, \ldots, n$, define $y^i = \iota ( \xi^i )$. Thus $y^1,
\ldots, y^n
\in {\cal R}^n$ and $$ A_{ij} = |y^i - y^j|^2 \ \qed$$

\vskip 18pt

{\bf Remark.} Of course, the fact that $ y^n = 0 $ by this
construction is of no import; we may take any translate of the $n$
vectors $y^1,
\ldots, y^n$ if we wish.

\vskip 18pt

{\bf 2.\ Applications}

\vskip 12pt
\hskip 18pt In this section we introduce a class of functions inducing AND
matrices and then use our characterization theorem 1.4 to prove a
simple, but rather useful, theorem on composition within this class.
We illustrate these ideas in examples 2.3-2.5. The remainder of the
section then uses theorems 1.4 and 2.2 to deduce results concerning
powers of the Euclidean norm. This enables us to derive the promised
$p$-norm result in theorem 2.11.

\proclaim {Definition 2.1}. We shall call $f:[0,\infty) \rightarrow [0,\infty) $ a {\bf conditionally
negative definite function of {order 1} (CND1)} if, for any positive
integers $n$ and $d$, and for any points $ x^1, \ldots, x^n
\in {\cal R}^d$, the matrix $A \in {\cal R}^{n \times n} $
defined by $$ A_{ij} = f(|x^i - x^j|^2), \hbox{ for } 1 \le i, j \le
n, $$ is AND. Furthermore, we shall call $f$ {\bf strictly CND1} if
the matrix $A$ is strictly AND whenever $n \ge 2$ and the points $x^1,
\ldots, x^n$ are distinct.

This terminology follows that of Micchelli [4], definition 2.1 .  We
see that the matrix $A$ of the previous definition satisfies the
conditions of proposition 1.3 if $f$ is strictly CND1, $n \ge 2$ and
the points $x^1, \ldots, x^n$ are distinct.

\vskip 18pt

{\bf Theorem 2.2.}

{ \sl (1) Suppose that $f$ and $g$ are CND1 functions and that $f(0) =
0 $. Then $g \circ f$ is also a CND1 function. Indeed, if $g$ is
strictly CND1 and f vanishes only at $0$, then $g \circ f$ is strictly
CND1. }

{\sl (2) Let A be an AND matrix with all diagonal entries zero.  Let
$g$ be a CND1 function.  Then the matrix defined by $$ B_{ij} =
g(A_{ij}), \hbox{ for } 1 \le i, j \le n,$$ is AND. Moreover, if $n
\ge 2$ and no off-diagonal elements of $A$ vanish, then $B$ is
strictly AND whenever $g$ is strictly AN. }

\vskip 12pt

{\it Proof.}

(1)The matrix $ A_{ij} = f(|x^i - x^j|^2)$ is an AND matrix with all
diagonal entries zero.  Hence, by theorem 1.4, we can find $n$ vectors
$ y^1, \ldots, y^n \in {\cal R}^n $ such that $$ f(|x^i - x^j|^2) =
|y^i - y^j|^2. $$

But g is a CND1 function, and so the matrix $B \in {\cal R}^{n
\times n} $ defined by $$ B_{ij} = g(|y^i - y^j|^2) = g \circ f(|x^i -
x^j|^2),$$ is also an AND matrix. Thus $g \circ f$ is a CND1 function.

The condition that $f$ vanishes only at $0$ allows us to deduce that
$y^i \ne y^j$, whenever $i \ne j$. Thus $B$ is strictly AND if $g$ is
strictly CND1.

(2) We observe that $A$ satisfies the hypotheses of theorem 1.4. We
may therefore write $A_{ij} = |y^i - y^j|^2$, and thus $B$ is AND
because $g$ is CND1. Now, if $A_{ij} \ne 0$ if $i \ne j$, then the
vectors $y^1, ...  , y^n$ are distinct, so that $B$ is strictly AND if
$g$ is strictly CND1 \qed

\vskip 18pt

\hskip 18pt For the next two examples only, we shall need the following concepts.
Let us call a function $g:[0,\infty)
\rightarrow [0,\infty) $ {\bf positive definite} if, for any positive integers
$n$ and $d$, and for any points $ x^1, \ldots, x^n \in {\cal R}^d$,
the matrix $A \in {\cal R}^{n \times n} $ defined by $$ A_{ij} =
g(|x^i - x^j|^2), \hbox{ for } 1 \le i, j \le n, $$ is non-negative
definite.  Furthermore, we shall call $g$ {\bf strictly positive
definite} if the matrix $A$ is positive definite whenever the points $
x^1, \ldots, x^n $ are distinct.  We reiterate that these last two
definitions are needed only for examples 2.3 and 2.4.

\vskip 12pt

{\bf Example 2.3.} A Euclidean distance matrix $A$ is AND, indeed
strictly so given distinct points. This was proved by Schoenberg [7]
and rediscovered by Micchelli [4]. Schoenberg also proved the stronger
result that the matrix
$$ A_{ij} = |x^i - x^j|^\alpha, \hbox{ for } 1 \le i,j \le n,$$
is strictly AND given distinct points $x^1, \ldots, x^n \in {\cal R}^d$, $n \ge 2$
and $0 < \alpha < 2$. We shall derive this fact using Micchelli's methods in
corollary 2.7 below, but we shall use the result here to
illustrate theorem 2.2. 
We see that, by theorem 1.4, there exist $n$ vectors $ y^1, \ldots, y^n
\in {\cal R}^n $ such that $$ A_{ij} \equiv |x^i - x^j|^\alpha = |y^i -
y^j|^2. $$

The vectors $ y^1, \ldots, y^n$ must be distinct whenever the points $
x^1, \ldots, x^n \in {\cal R}^d$ are distinct, since $A_{ij} \ne 0$
whenever $i \ne j$.

\hskip 12pt Now let $g$ denote any strictly positive definite function. Define $B \in
{\cal R}^{n \times n} $ by $$ B_{ij} \equiv g(A_{ij}). $$ Thus $$
g(|x^i - x^j|^\alpha) = g(|y^i - y^j|^2) .$$ Since we have shown that the
vectors $ y^1, \ldots, y^n$ are distinct, the matrix $B$ is therefore
positive definite.

\hskip 12pt For example, the function $g(t) = \exp(-t)$ is a strictly positive
definite function.  For an elementary proof of this fact, see Micchelli
[4], p.15 .  Thus the
matrix whose elements are $$ B_{ij} = \exp(\ -|x^i - x^j|^\alpha), 1 \le i,j
\le n, $$ is always (i) non-negative definite, and (ii) positive
definite whenever the points $ x^1, \ldots, x^n $ are distinct
\qed

\vskip 12pt

{\bf Example 2.4.} This will be our first example using a $p$-norm
with $p \ne 2$. Suppose we are given distinct points $ x^1, \ldots,
x^n
\in {\cal R}^d$. Let us define $A \in {\cal R}^{n \times n}
$ by $$ A_{ij} = \Vert x^i - x^j \Vert_1 . $$

Furthermore, for $k = 1, \ldots, d$, let $A^{(k)} \in {\cal R}^{n
\times n}$ be given by $$ A_{ij}^{(k)} = | x^i_k - x^j_k |,$$
recalling that $ x^i_k$ denotes the $k^{th}$ coordinate of the point
$x^i$.

\hskip 12pt We now remark that $ A = \sum_{i=1}^d A^{(k)} $. But every $A^{(k)}$
is a Euclidean distance matrix, and so every $A^{(k)}$ is AND.
Consequently $A$, being the sum of AND matrices, is itself AND.  Now
$A$ has all diagonal entries zero. Thus, by theorem 1.4, we can
construct $n$ vectors $ y^1, \ldots, y^n \in {\cal R}^n $ such that $$
A_{ij} \equiv \Vert x^i - x^j \Vert_1 = |y^i - y^j|^2. $$

As in the preceding example, whenever the points $ x^1, \ldots, x^n $
are distinct, so too are the vectors $ y^1, \ldots, y^n$.

\vskip 12pt
This does not mean that $A$ is non-singular.  Indeed, Dyn, Light and
Cheney [2] observe that the 1-norm distance matrix is singular for the
distinct points $\{ (0,0), (1,0), (1,1), (0,1) \} $.
\vskip 12pt
\hskip 12pt Now let $g$ be any strictly positive definite function. Define $B \in
{\cal R}^{n \times n} $ by $$ B_{ij} = g(A_{ij}) = g( \Vert x^i - x^j
\Vert_1 ) = g(|y^i - y^j|^2) .$$

Thus $B$ is positive definite.

\hskip 12pt For example, we see that the matrix $ B_{ij} = \exp(\ -\Vert x^i - x^j
\Vert_1 )$ is positive definite whenever the points $ x^1, \ldots, x^n $
are distinct \qed

\vskip 12pt

{\bf Example 2.5.} As in the last example, let $ A_{ij} = \Vert x^i -
x^j \Vert_1 $, where $n \ge 2$ and the points $ x^1, \ldots, x^n $ are
distinct.  Now the function $ f(t) = (1+t)^{1 \over 2}$ is strictly
CND1 (\ Micchelli [4]\ ).  This is the CND1 function generating the
multiquadric interpolation matrix.  We shall show the matrix $B \in
{\cal R}^{n \times n} $ defined by $$ B_{ij} = f(A_{ij}) = (1+\Vert
x^i - x^j \Vert_1)^{1 \over 2}$$ to be strictly AND.

\hskip 12pt Firstly, since the points $ x^1, \ldots, x^n $ are distinct, the previous
example shows that we may write $$ A_{ij} = \Vert x^i - x^j \Vert_1 =
|y^i - y^j|^2 ,$$ where the vectors $y^1, \ldots, y^n$ are distinct.
Thus, since $f$ is strictly CND1, we deduce from definition 2.1 that
$B$ is a strictly AND matrix \qed

\vskip 12pt
%** Insert new lemmata here replacing 2.6 etc

\hskip 12pt We now return to the mainstream of the paper.
Recall that a function $f$ is completely monotonic provided
that
$$ (-1)^k f^{(k)}(x) \ge 0 , \hbox{ for every }k = 0, 1, 2, \ldots
\hbox{ and for } 0 < x < \infty .$$
We now require a theorem of Micchelli [4], restated in our notation.
\vskip 12pt
\proclaim{Theorem 2.6}. Let $f:[0,\infty) \rightarrow [0,\infty)$ have a completely monotonic
derivative. Then $f$ is a CND1 function. Further, if $f^\prime$ is non-constant, then
$f$ is strictly CND1.

{\it Proof. } This is theorem 2.3 of Micchelli [4] \qed

\vskip 12pt

\proclaim{Corollary 2.7}. The function $g(t) = t^\tau$ is strictly CND1 for every $\tau \in (0,1)$.

{\it Proof. }The conditions of the previous theorem are satisfied by $g$ \qed

\vskip 12pt

We see now that we may use this choice of $g$ in theorem 2.2, as in the
following corollary.

\vskip 12pt

\proclaim{Corollary 2.8}. For every $\tau \in (0,1)$ and for every positive integer $k \in [1, d]$, define
$A^{(k)} \in {\cal R}^{n \times n} $ by $$ A_{ij}^{(k)} = |x^i_k -
x^j_k|^{2\tau }, \hbox{ for } 1 \le i,j \le n.$$ Then every $A^{(k)}$
is AND.

{\it Proof. } For each $k$, the matrix $(|x^i_k-x^j_k|)_{i,j=1}^n $
is a Euclidean distance matrix. 
Using the function $g(t) = t^\tau$, we
now apply theorem 2.2 (2) to deduce that $A^{(k)} = 
g( |x^i-x^j|^2)$ is AND \qed
\vskip 12pt
We shall still use the notation $\Vert . \Vert_p$ when $p
\in (0,1)$, although of course these functions are not norms .
\vskip 12pt

{\bf Lemma 2.9. } {\sl For every $p \in (0,2)$, the matrix $A
\in {\cal R}^{n \times n} $ defined by $$ A_{ij} = \Vert x^i -
x^j \Vert_p^p, \hbox{ for } 1 \le i, j \le n, $$ is AND.  If $n \ge 2
$ and the points $x^1, \ldots, x^n$ are distinct, then we can find
distinct $y^1, \ldots, y^n \in {\cal R}^n$ such that $$ \Vert x^i -
x^j \Vert_p^p = |y^i - y^j|^2 . $$ }

{\it Proof.} If we set $p = 2\tau $, then we see that $\tau \in (0,1)$
and $A = \sum_{k=1}^d A^{(k)}$, where the $A^{(k)}$ are those matrices
defined in corollary 2.8. Hence so that each $A^{(k)}$ is AND, and
hence so is their sum.  Thus, by theorem 1.4, we may write $$A_{ij} =
\Vert x^i - x^j
\Vert_p^p = |y^i - y^j|^2 . $$ Furthermore, if $n \ge 2 $ and the points
$x^1, \ldots, x^n$ are distinct, then $A_{ij} \ne 0$ whenever $i \ne
j$, so that the vectors $y^1, \ldots, y^n$ are distinct \qed

\vskip 12pt
\proclaim {Corollary 2.10}. For any $p \in (0,2)$ and for any $\sigma \in (0,1)$, define 
$B \in {\cal R}^{n \times n} $ by $$B_{ij} = ( \Vert x^i - x^j
\Vert_p^p )^\sigma .$$ Then $B$ is AND. As before, if $n \ge 2 $ and
the points $x^1, \ldots, x^n$ are distinct, then $B$ is strictly AND.

{\it Proof.} Let $A$ be the matrix of the previous lemma and 
let $g(t) = t^\tau$. We now apply theorem 2.2 (2)
\qed
\vskip 12pt

\proclaim {Theorem 2.11}. For every $p \ \in (1,2)$, the $p$-norm distance matrix $B \in {\cal R}^{n \times n} $,
that is: $$B_{ij} = \Vert x^i - x^j \Vert_p , \hbox{ for } 1 \le i, j
\le n,$$ is AND. Moreover, it is strictly AND if $n
\ge 2 $ and the points $x^1, \ldots, x^n$ are distinct, in which case $$ (-1)^{n-1}\det B > 0 .$$

{\it Proof. } If $p \in (1,2)$, then $\sigma \equiv 1/p \ \in (0,1)$.
Thus we may apply corollary 2.12. The final inequality follows from
the statement of proposition 1.3 \qed
\vskip 12pt
\hskip 12pt We may also apply theorem 2.2 to the $p-$norm distance matrix, for $p \in (1,2]$, or indeed to the $p^{th}$ power
of the $p-$norm distance matrix, for $p \in (0,2)$. Of course, we do not have a norm
for $0<p<1$, but we define the function in the obvious way.
We need only note
that, in these cases, both classes satisfy the conditions of theorem
2.2 (2). We now state this formally for the $p-$norm distance matrix
\vskip 12pt
\proclaim {Corollary 2.12}. Suppose the matrix $B$ is the $p-$norm distance matrix defined in theorem 2.13. Then,
if $g$ is a CND1 function, the matrix $g(B)$ defined by $$ g(B)_{ij}
= g(B_{ij}), \hbox{ for } 1 \le i,j \le n,$$ is AND. Further, if $n
\ge 2$ and the points $x^1, \ldots, x^n$ are distinct, then $g(B)$ is
strictly AND whenever $g$ is strictly AN.

{\it Proof.} This is immediate from theorem 2.11 and the statement of
theorem 2.2 (2) \qed

\vskip 18pt

{\bf 3. The Case ${\bf p > 2}$}

\vskip 6pt

\hskip 12pt We are unable to use the ideas developed in the previous section to
understand this case. However, numerical experiment suggested the
geometry described below, which proved surprisingly fruitful.  We
shall view ${\cal R}^{m+n}$ as two orthogonal slices ${\cal R}^m
\oplus {\cal R}^n$. Given any $p > 2$, we take the vertices $\Gamma_m$
of $[-m^{-1/p}, m^{-1/p}]^m \subset {\cal R}^m$ and embed this in
${\cal R}^{m+n}$. Similarly, we take the vertices $\Gamma_n$ of
$[-n^{-1/p}, n^{-1/p}]^n \subset {\cal R}^n$ and embed this too in
${\cal R}^{m+n}$. We see that we have constructed two orthogonal cubes
lying in the $p$-norm unit sphere.

\vskip 12pt
{\bf Example.} If $m = 2$ and $n = 3$, then $\Gamma_m = \{ (\pm \alpha
,
\pm \alpha , 0, 0, 0) \}$ and $\Gamma_n = \{ (0, 0, \pm \beta, \pm \beta , \pm \beta ) \}
$, where $\alpha = 2^{-1/p}$ and $\beta = 3^{-1/p}$.

\vskip 12pt

\hskip 18pt Of course, given $m$ and $n$, we are interested in values of $p$ for which the 
$p-$norm distance matrix generated by $\Gamma_m \cup \Gamma_n$ is
singular. Thus we ask whether there exist scalars $\{\lambda_y\}_{\{y
\in \Gamma_m\}}$ and $\{\mu_z\}_{\{z \in \Gamma_n\}}$, not all zero,
such that the function $$ s(x) = \sum_{y \in \Gamma_m} \lambda_y \Vert
x-y \Vert_p + \sum_{z \in \Gamma_n} \mu_z \Vert x-z \Vert_p $$
vanishes at every interpolation point. In fact, we shall show that
there exist scalars $\lambda$ and $\mu$, not both zero, for which the
function $$ s(x) = \lambda \sum_{y \in \Gamma_m} \Vert x-y \Vert_p +
\mu \sum_{z \in \Gamma_n} \Vert x-z \Vert_p $$ vanishes at every
interpolation point.

\hskip 12pt We notice that

\hskip 12pt (i) \ For every $y \in \Gamma_m$ and $z \in \Gamma_n$,
we have $\Vert y - z \Vert_p = 2^{1/p}$.

\hskip 12pt (ii) The sum $\sum_{y \in \Gamma_m} \Vert \tilde y - y \Vert_p$
takes the same value for every vertex $\tilde y \in \Gamma_m$, and
similarly, {\it mutatis mutandis}, for $\Gamma_n$.

\hskip 12pt Thus our interpolation equations reduce to two in number: $$ \lambda
\sum_{y \in \Gamma_m} \Vert \tilde y - y \Vert_p
\ + \ 2^{n+1/p}\mu \ = \ 0,$$
and $$ 2^{m+1/p}\lambda \ + \ \mu \sum_{z \in \Gamma_n} \Vert
\tilde z - z \Vert_p \ = \ 0 , $$ where by (ii) above, we see that
$\tilde y $ and $\tilde z$ may be any vertices of $\Gamma_m, \Gamma_n$
respectively.

\hskip 12pt We now simplify the (1,1) and (2,2) elements of our reduced system by
use of the following lemma.
\vskip 12pt
{\bf Lemma 3.1. }{\sl Let $\Gamma $ denote the vertices of $[0, 1]^k
$.  Then } $$ \sum_{x \in
\Gamma } \Vert x \Vert_p = \sum_{l=0}^k {k \choose l} l^{1/p} .$$
\vskip 12pt
{\it Proof. } Every vertex of $\Gamma $ has coordinates taking the
values $0$ or $1$. Thus the distinct $p$-norms occur when exactly $l$
of the coordinates take the value $1$, for $l = 0, \ldots, k$; each of
these occurs with frequency $k \choose l$ \qed
\vskip 12pt 
{\bf Corollary 3.2.} {\sl $$ \sum_{y \in \Gamma_m} \Vert \tilde y - y
\Vert_p = 2 \sum_{k=0}^m {m \choose k} (k/m)^{1/p}, \hbox{ for every }
\tilde y \in \Gamma_m, \hbox{ and }$$
$$ \sum_{z \in \Gamma_n} \Vert \tilde z - z \Vert_p = 2 \sum_{l = 0}^n
{n \choose l} (l/n)^{1/p}, \hbox{ for every }
\tilde z \in \Gamma_n.$$ }
\vskip 12pt
{\it Proof. } We simply scale the result of the previous lemma by
$2m^{-1/p}$ and $2n^{-1/p}$ respectively \qed
\vskip 12pt

\hskip 12pt With this simplification, the matrix of our system
becomes $$ \left(\matrix{ 2\sum_{k=0}^m {m \choose k} (k/m)^{1/p} & \
&2^n.2^{1/p}\cr
\ & \ & \ \cr          
2^m.2^{1/p}& \ & 2\sum_{l=0}^n {n \choose l} (l/n)^{1/p}\cr}
\right).$$

\hskip 12pt We now recall that
$$ B_i(f_p, 1/2) = 2^{-i} \sum_{j=0}^i {i \choose j} (j/i)^{1/p}$$ is
the Bernstein polynomial approximation of order $i$ to the function
$f_p(t) = t^{1/p}$ at $t = 1/2$. Our reference for properties for
Bernstein polynomial approximation will be Davis [1], sections 6.2 and
6.3. Hence, scaling the determinant of our matrix by $2^{-(m+n)}$, we
obtain the function $$\varphi_{m,n}(p) = 4\/B_m(f_p, 1/2)\/B_n(f_p,
1/2) - 2^{2/p} .$$ We observe that our task reduces to investigation
of the zeros of $\varphi_{m,n}$.

\hskip 12pt We first deal with the case $m = n$, noting the factorization:
$$ \varphi_{n,n}(p) = \{ 2B_n(f_p,1/2) + 2^{1/p}\} \{ 2B_n(f_p,1/2) -
2^{1/p}\} .$$ Since $f_p(t) \ge 0$, for $t \ge 0$ we deduce from the
monotonicity of the Bernstein approximation operator that $B_n(f_p,
1/2) \ge 0$. Thus the zeros of $\varphi_{n,n}$ are those of the factor
$$\psi_n(p) = 2B_n(f_p, 1/2) - 2^{1/p}.$$

\vskip 12pt
{\bf Proposition 3.3.} {\sl $\psi_n$ enjoys the following properties.}

\hskip 12pt{\sl (1) $\psi_n(p) \rightarrow \psi (p)$, where $\psi (p) =
2^{1 - 1/p} - 2^{1/p}$, as $n \rightarrow \infty$.

\hskip 12pt(2) For every $p > 1$, $\psi_n (p) < \psi_{n+1}(p)$, for every
positive integer $n$.

\hskip 12pt(3) For each $n$, $\psi_n$ is strictly increasing for $p \in 
[1, \infty)$.

\hskip 12pt(4) For every positive integer $n$, $\lim_{p \to \infty} \psi_n
(p) = 1 - 2^{1-n}$. }
\vskip 12pt
{\it Proof. }

\hskip 12pt(1) This is a consequence of the convergence of Bernstein polynomial
approximation.

\hskip 12pt(2) It suffices to show that $B_n(f_p, 1/2) < B_{n+1}(f_p, 1/2)$, for
$p > 1$ and $n$ a positive integer. We shall use Davis [1], theorem
6.3.4: If $g$ is a convex function on $[0,1]$, then $B_n(g,x) \ge
B_{n+1}(g,x)$, for every $x \in [0,1]$. Further, if $g$ is non-linear
in each of the intervals $[{ {j-1}\over n}, {j \over n}]$, for $j = 1,
\ldots, n$, then the inequality is strict.

\hskip 12pt Every function $f_p$ is concave and non-linear on $[0,1]$
for $p > 1$, so that this inequality is strict and reversed.

\hskip 12pt(3) We recall that 
$$ \psi_n(p) = 2B_n(f_p, 1/2) - 2^{1/p} = 2^{1-n} \sum_{k=0}^n {n
\choose k} (k/n)^{1/p} - 2^{1/p}.$$ Now, for $p_2 > p_1 \ge 1$, we
note that $t^{1/p_2} > t^{1/p_1}$, for $t \in (0,1)$, and also that
$2^{1/p_2} < 2^{1/p_1}$. Thus $(k/n)^{1/p_2} > (k/n)^{1/p_1}$, for $k
= 1, \ldots, n-1$ and so $\psi_n(p_2) > \psi_n(p_1)$.

\hskip 12pt(4) We observe that, as $p \rightarrow \infty$,
$$ \psi_n(p) \rightarrow 2^{1-n} \sum_{k=1}^n {n \choose k} - 1 = 2( 1
- 2^{-n}) - 1 = 1 - 2^{1-n} \qed $$

\vskip 12pt
{\bf Corollary 3.4.} {\sl For every integer $n > 1$, each $\psi_n$ has
a unique root $p_n \in (2, \infty )$.  Further, $p_n \rightarrow 2$
strictly monotonically as $n \rightarrow \infty$.}
\vskip 12pt
{\it Proof. } We first note that $\psi (2) = 0$, and that this is the
only root of $\psi$.  By proposition 3.3 (1) and (2), we see that
$$\lim_{n \rightarrow \infty}{\psi_n(2)} = \psi (2) = 0 \hbox{ and }
\psi_n(2) < \psi_{n+1}(2) < \psi (2) = 0.$$ By proposition 3.3 (4), we
know that, for $n > 1$, $\psi_n$ is positive for all sufficiently
large $p$.  Since every $\psi_n$ is strictly increasing by proposition
3.3 (3), we deduce that each $\psi_n$ has a unique root $p_n \in (2,
\infty)$ and that $\psi_n(p) < (>) 0 $ for $p < (>) p_n$.

\hskip 12pt We now observe that $\psi_{n+1}(p_n) > \psi_n(p_n) = 0$, by proposition 3.3 (2), whence
$2 < p_{n+1} < p_n$.  Thus $(p_n)$ is a monotonic decreasing sequence
bounded below by $2$.  Therefore it is convergent with limit in $[2,
\infty)$.  Let $p^*$ denote this limit.  To prove that $p^* = 2$, it
suffices to show that $\psi(p^*) = 0$, since $2$ is the unique root of
$\psi$.  
Now suppose that $\psi(p^*) \ne 0$. By continuity, $\psi$ is bounded away from zero
in some compact neighbourhood $N$ of $p^*$. We now recall the following theorem
of Dini: If we have a monotonic increasing sequence of continuous real-valued functions
on a compact metric space with continuous limit function, then the convergence is uniform.
A proof of this result may be found in many texts, for example Hille [3], p. 78.
Thus $\psi_n \rightarrow \psi$ uniformly in $N$. Hence there is an integer $n_0$ 
such that $\psi_n$ is bounded away from zero for every $n \ge n_0$. But 
$p^* = \lim p_n$ and $\psi_n(p_n) = 0$ for each $n$, so that we have reached
a contradiction. Therefore $\psi(p^*) = 0$ as required \qed

\vskip 12pt
\hskip 12pt Returning to our original scaled determinant $\varphi_{n,n}$, we see that
$\Gamma_n \cup \Gamma_n$ generates a singular $p_n$-norm distance
matrix and $p_n \searrow 2$ as $n \rightarrow \infty$.  Furthermore $$
\varphi_{m,m}(p) < \varphi_{m,n}(p) < \varphi_{n,n}(p) , \hbox{ for }
1 < m < n, $$ using the same method of proof as in proposition 3.3
(2).  Thus $\varphi_{m,n}$ has a unique root $p_{m,n}$ lying in the
interval $(p_n, p_m).$ We have therefore proved the following theorem.
\vskip 12pt
{\bf Theorem 3.5. }{\sl For any positive integers $m$ and $n$, both
greater than $1$, there is a $p_{m,n} > 2$ such that the $\Gamma_m
\cup
\Gamma_n$-generated $p_{m,n}$-norm distance matrix is singular.
Furthermore, if $1 < m < n$, then $$ p_m \equiv p_{m,m} > p_{m,n} >
p_{n,n} \equiv p_n, $$ and $p_n \searrow 2 \hbox{ as } n \rightarrow
\infty$.}
\vskip 12pt
\hskip 12pt Finally, we deal with the ``gaps'' in the sequence $(p_n)$ as follows.
Given a positive integer $n$, we take the configuration $\Gamma_n \cup
\Gamma_n (\vartheta )$, where $\Gamma_n (\vartheta )$ denotes the vertices
of the scaled cube $[-\vartheta n^{-1/p}, \vartheta n^{-1/p}]^n$ and
$\vartheta > 0$. The $2 \times 2$ matrix deduced from corollary 3.2 on
page 8 becomes $$ \left(\matrix{ 2 \sum_{k=0}^n {n \choose k}
(k/n)^{1/p} & \ &2^n (1+\vartheta^p)^{1/p}\cr
\ & \ & \ \cr
2^n (1+\vartheta^p)^{1/p}& \ & 2 \vartheta \sum_{k=0}^n {n \choose k}
(k/n)^{1/p}\cr}
\right).$$
Thus, instead of the function $\varphi_{n,n}$ discussed above, we now
consider its analogue: $$ \varphi_{n,n,\vartheta}(p) = 4\vartheta
B_n^2(f_p, 1/2) - (1 + \vartheta^p)^{2/p}. $$ If $p > p_n$, the unique
zero of our original function $\varphi_{n,n}$, we see that
$\varphi_{n,n,1}(p)
\equiv \varphi_{n,n} (p) > 0 $, because every $\varphi_{n,n}$ is strictly
increasing, by proposition 3.3 (3). However, we notice that
$\lim_{\vartheta \to 0} \varphi_{n,n,\vartheta}(p) = -1$, so that
$\varphi_{n,n,\vartheta}(p) < 0$ for all sufficiently small $\vartheta
> 0$.  Thus there exists a $\vartheta^* > 0 \hbox{ such that }
\varphi_{n,n,\vartheta^* }(p) = 0$. Since this is true for every $p >
p_n$, we have strengthened the previous theorem. We now state this
formally.
\vskip 12pt
{\bf Theorem 3.6. }{\sl For every $p > 2$, there is a configuration of
distinct points generating a singular $p$-norm distance matrix. }

\vskip 12pt
\hskip 12pt It is interesting to investigate how rapidly the sequence of zeros
$(p_n)$ converges to $2$. We shall use Davis [1], theorem 6.3.6, which
states that, for any bounded function $f$ on $[0,1]$,
$$\lim_{n\to\infty}{n(B_n(f,x) - f(x)) } = {1 \over 2} x(1-x)f^{\prime
\prime}(x),
\hbox{ whenever }f^{\prime \prime }(x) \hbox{ exists}. $$Applying this to 
$$\psi_n(p) = 2 B_n(f_p, 1/2) - 2^{1/p}, $$ we shall derive the
following bound.
\vskip 12pt
{\bf Proposition 3.7. }$p_n = 2 + O(n^{-1}).$
\vskip 6pt
{\it Proof.} We simply note that $$\eqalign { 0 &= \psi_n(p_n) \cr &=
\psi (p_n) + O(n^{-1}), \hbox{ by Davis [1] 6.3.6,} \cr &=
\psi (2) + (p_n-2) \psi^\prime (2) + {\it o}(p_n-2) + O(n^{-1}) .}$$
Since $\psi^\prime (2) \ne 0$, we have $p_n - 2 = O(n^{-1})$ \qed

\vskip 18pt
{\bf 4. Acknowledgments}
\vskip 6pt
\hskip 12pt The author is indebted to his research supervisors
in Cambridge and Harwell, namely Prof. M. J. D. Powell and Dr. N. I.
M. Gould, for their numerous helpful criticisms and encouragement. In
particular, the author is grateful to M. J. D. Powell for his
suggestions concerning the proof of theorem 3.5. Thanks are
also due to I. J. Leary for several helpful conversations and to
the referee for many helpful and necessary suggestions. It was the
referee who brought the paper of Schoenberg [7] to my attention.

\vskip 18pt
{\bf 5. References}

\vskip 6pt
[1] P. J. Davis (1975): Interpolation and Approximation. Dover
Publications, New York.
\vskip 12pt
[2] N. Dyn, W. A. Light and E. W. Cheney (1989): {\sl Interpolation by
piecewise linear radial basis functions.} To appear in J. Approx.
Theory.
\vskip 12pt
[3] E. Hille (1962): Analytic Function Theory, Vol. II. Waltham, Massachusetts:
Ginn and Co.
\vskip 12pt
[4] C. A. Micchelli (1986): {\sl Interpolation of scattered data:
distance matrices and conditionally positive functions.} Constructive
Approximation {\bf 2}:11-22.
\vskip 12pt
[5] I. J. Schoenberg (1935): {\sl Remarks to Maurice Fr\'echet's
article ``Sur la definition d'une classe d'espace distanci\'es
vectoriellement applicable sur l'espace d'Hilbert.''} Ann. of
Math., {\bf 36}: 724-732.
\vskip 12pt
[6] I. J. Schoenberg (1937): {\sl On certain metric spaces arising
from Euclidean space by a change of metric and their embedding in
Hilbert space.} Ann. of Math., {\bf 38}: 787-793.
\vskip 12pt
[7] I. J. Schoenberg (1938): {\sl Metric spaces and completely
monotone functions.} Ann. of Math., {\bf 39}: 811-841.

\bye